\documentclass[12pt,reqno]{amsart}

\title[Multipliers]{Multipliers on a Hilbert space of functions on $\R$}

\author[Violeta Petkova]{Violeta Petkova}

\address{LMAM, Universit\'e de Metz UMR 7122,Ile du Saulcy 57045 Metz Cedex 1, France.}

\email{petkova@univ-metz.fr}

\usepackage{amsmath} 

\usepackage{amssymb}

\def\squarebox#1{\hbox to #1{\hfill\vbox to #1{\vfill}}}

\newcommand{\F}{{\mathcal F}}

\newcommand{\I}{[-\ln \rho (S^{-1}), \:\ln \rho(S)]}

\newcommand{\M}{{\mathcal M}}

\newcommand{\R}{{\mathbb R}}

\newcommand{\LL}{L_{loc}^1(\R)}

\newcommand{\C}{{\mathbb C}}

\newcommand{\N}{{\mathbb N}}

\newcommand{\w}{{\omega}}

\newcommand{\lw}{{L_\omega^2({\R})}}

\newcommand{\Cc}{{C_c(\R)}}

\renewcommand{\Re}{\mathop{\rm Re}\nolimits}

\renewcommand{\Im}{\mathop{\rm Im}\nolimits}

\theoremstyle{plain}

\newtheorem{thm}{Theorem}

\newtheorem{deff}{Definition}

\theoremstyle{definition}

\vspace{1cm}

\numberwithin{equation}{section}

\begin{document}

%\address{Violeta Petkova\\

%LMAM, \\

%Metz, France.\\}

%\email{petkova@univ-metz.fr}

\maketitle

\begin{abstract} For a Hilbert space $H \subset L^1_{loc} (\R)$ of functions on $\R$ we obtain a representation theorem for the multipliers $M$ commuting with the shift operator $S$. This generalizes the classical result
for multipliers in $L^2(\R)$ as well as our previous result for multipliers in weighted space $L^2_{\omega}(\R).$ Moreover, we obtain a description of the spectrum of $S$. 

\end{abstract}


\vspace{0.6cm}
\section{Introduction}

Let $H\subset L^1_{loc}(\R)$ be a Hilbert space of functions on $\R$ with values in $\C$. Denote by $\|.\|$ (resp. $<.,.>$) the norm (resp. the scalar product) on $H$.
Let $\Cc$ be the set of continuous functions on $\R$ with compact support. 
For a compact $K$ of $\R$ denote by $C_K(\R)$ the subset of functions of $\Cc$ with support in $K$ and denote by $\hat{f}$ or by $\F(f)$ the usual Fourier transform of $f\in L^2(\R)$.
Let $S_x$ be the operator of translation by $x$ defined on $H$ by 
$$(S_xf)(t)=f(t-x),\:a.e.\:\: t \in \R.$$
Let $S$ (resp. $S^{-1}$) be the translation by 1 (resp. -1). Introduce the set 
$$\Omega=\Big\{z\in \C,\: -\ln \rho(S^{-1})\leq \Im z\leq \ln \rho(S)\Big\},$$
where $\rho(A)$ is the spectral radius of $A$ and
let ${I}$ be the interval $\I$.
Assuming the identity map $i:H \longrightarrow \LL$ continuous,
 it follows from the closed graph theorem that if $S_x(H)\subset H$, for $x \in \R$, then the operator $S_x$ is bounded from $H$ into $H$. In this paper we suppose that $H$ satisfies the following conditions:\\

 %Likewise, if $\Gamma_\chi(E)\subset E$, for $\chi \in \G$, the operator $ \Gamma_\chi$ is bounded from $E$ into $E$. We will be interested in Banach spaces $E$ satisfying the following conditions:\\

(H1) $\Cc \subset H \subset \LL$, with continuous inclusions, and $\Cc$ is dense in $H$.\\

(H2) For every $x \in \R$, $S_x(H)\subset H$ and $\sup_{x\in K}\|S_x\|<+\infty$, for every compact set $K\subset \R$. \\

(H3) For every $\alpha \in \R$ let $T_\alpha$ be the operator defined by 
$$T_\alpha:H\ni f(x) \longrightarrow f(x)e^{i\alpha x}, \: x \in \R.$$
 We have $T_\alpha(H)\subset H$ and, moreover,
 $\sup_{\alpha \in \R}\| T_\alpha \|<+\infty.$\\
 
 (H4) There exists $C > 0$ and $a \geq 0$ such that $\|S_x\|\leq Ce^{a|x|},\:\forall x\in\R.$\\
 
 Set $|||f|||=\sup_{\alpha \in \R} \|T_\alpha f\|,$ for $f \in H$. The norm $|||.|||$ is equivalent to the norm of $H$ and without loss of generality, we can consider below that $T_\alpha$ is an isometry on $H$ for every $\alpha\in \R$.
Obviously, the condition $(H3)$ holds for a very large class of Hilbert spaces.\\
We give some examples of Hilbert spaces satisfying our hypothesis.\\
{\bf Example 1.} A weight $\w$ on $\R$ is a non negative function on $\R$ such that
$$\sup_{x\in \R} \frac{\w(x+y)}{\w(x)}<+\infty, \:\forall y\in \R.$$
Denote by $\lw$ the space of measurable functions on $\R$ such that 
$$\int_\R |f(x)|^2\w(x)^2dx<+\infty.$$
The space $\lw$ equipped with the norm
$$\|f\|=\Big(\int_\R |f(x)|^2 \w(x)^2 dx\Big)^{\frac{1}{2}}$$
is a Hilbert space satisfying our conditions (H1)-(H3). Moreover, we have the estimate
\begin{equation} \label{eq:sg}
\|S_{t}\|\leq Ce^{m |t|},\:\forall t\in \R,
\end{equation}
where $C > 0$ and $m \geq 0$ are constants. 
This follows from the fact that $\w$ is equivalent to the special weight $\w_0$ constructed in \cite{M}. The details of the construction of $\w_0$ are given in \cite{V1}, \cite{M}. Below after Theorem 2 we give some examples of weights. \\

%{\bf Example 2.}
%Let $A$ be a real continuous function on $[0, +\infty[$, such that $A(0)=0$ and assume that the function
%$$\R\ni x\longrightarrow \frac{A(x)}{x}$$
% is increasing on $\R^+$. 
%Let $L_A(\R)$ be the spaces of measurable functions on $\R$ such that
%$$\int_\R A\Big( \frac{|f(x)|}{t}\Big)dx<+\infty,$$
%for a $t>0$ and let
%$$\|f\|_A=\inf\Big\{t>0\:|\:\int_\R A\Big( \frac{|f(x)|}{t}\Big)dx\leq 1\Big\}$$
%for $f\in L_A(\R)$. Then $L_A(\R)$ is a Banach space called a Birnbaum-Orlicz space (cf. \cite{O}). For every $t\in \R$, we have $\|S_{t}\|=1$. 
%It is clear that $L_A(\R)$ satisfies (H1), (H2), (H3) and (H4). 
%For some functions $A$, $L_A(\R)$ is a Hilbert space. \\

\begin{deff}
A bounded operator $M$ on $H$ is called a multiplier if
$$MS_x=S_xM,\:\forall x\in \R.$$
\end{deff}
Denote by $\M$ the algebra of the multipliers. 
Our aim is to obtain a representation theorem for multipliers on $H$ and to characterize the spectrum of $S$. 
These two problems are closely related. In \cite{V1} we have obtained a representation theorem for multipliers on $\lw$. Here we generalize our result for multipliers on a Hilbert space and shift operators satisfying the conditions $(H1)-(H4)$. Our proof is shorter than that in \cite{V1}. The main improvement is based on an application of the link between the spectrum $\sigma(S_t)$ of a element of the group $(S_t)_{t\in \R}$ and the spectrum $\sigma(A)$ of the generator $A$ of this group. In general, in the setup we deal with the spectral mapping theorem 
$$\sigma(S_t) \setminus \{0\} = e^{\sigma(tA)}$$ 
is not true. To establish the crucial estimate in Theorem 4 we use the general results (see \cite{G} and \cite{H})
for the characterization of the spectrum of $S_t$ by the behavior of the resolvent of $A.$ This idea has been used in
\cite{V7} for $L_{\omega}^2(\R)$ but one point in our argument needs a more precise proof and in this paper we do this in the general case.\\

 Denote by $(f)_a$ the function
$$\R\ni x\longrightarrow f(x)e^{ax}.$$
We prove the following
\begin{thm}

For every $M\in \M$, and for every $a\in I= [-\ln \rho(S^{-1}),\:\ln \rho(S)],$ we have\\
1) $(Mf)_a\in L^2(\R),\:\forall f\in \Cc$.\\
2) There exists $\mu_{(a)}\in L^\infty(\R)$ such that 
$$\int_\R (Mf)(x)e^{ax}e^{-itx}dx=\mu_{(a)}(t)\int_\R f(x)e^{ax}e^{-itx}dx,\:a.e. $$
i.e.\\
$$\widehat{(Mf)_a}=\mu_{(a)} \widehat{(f)_a}.$$
3) If $\overset{\circ}{I} \neq \emptyset $ then the function $\mu(z)=\mu_{(\Im z)}(\Re z)$ is holomorphic on $\overset{\circ}\Omega$.

\end{thm}

\begin{deff}
Given $M\in \M$, if $\overset{\circ}{\Omega}\neq \emptyset$, we call symbol of $M$ the function $\mu$ defined by 
$$\mu(z)=\mu_{(\Im z)}(\Re z),\:\forall z\in \overset{\circ}{\Omega}.$$
 Moreover, if $a=-\ln \rho(S^{-1})$ or $a=\ln \rho(S)$, the symbol $\mu$ is defined for $z=x+ia$ by the same formula for almost all $x\in \R$.
\end{deff}
Denote by $\sigma(A)$ the spectrum of the operator $A$. 
From Theorem 1 we deduce the following interesting spectral result.

\begin{thm}
We have
$$\sigma(S)=\Big\{z\in \C:\:\frac{1}{\rho(S^{-1})}\leq |z|\leq \rho(S)\Big\}.$$
\end{thm}
To prove this characterization of the spectrum of $S$ we exploit the existence of a symbol of every multiplier. Notice that in general $S$ is not a normal operator and there are no spectral calculus which could characterize the spectrum of $S$. On the other hand, Theorem 2 has been used in \cite{V8} to obtain spectral mapping theorems for a class of multipliers.
Now we give some examples of weights.\\

{\bf Example 2.} The function $ \w(x)=e^x$ is a weight. For the associated weighted space $\lw$ we obtain $\sigma(S)=\{z\in\C,\: |z|=e\}.$\\

{\bf Example 3.} The functions of the form $\w(x)=1+|x|^\alpha$, for $\alpha \in \R$ are weights and we get $\sigma(S)=\{z\in \C,\: |z|=1\}$.\\

{\bf Example 4.} Let $\w(x)=e^{a|x|^b}$ with $a>0$ and $0<b<1$. Then in $\lw$ we have 
$$\sigma(S)= \{z\in \C,\: e^{-a}\leq |z|\leq e^a\}.$$

{\bf Example 5.} Functions like 
$$e^{\frac{|x|}{\ln (2+|x|)}},\:\:\: e^{|x|}(1+|x|^2)^n,\:{\rm for}\: n>0$$ 
also are weights.\\ 
The weights in the Examples 4 and 5 are used to illustrate Beurling algebra theory (cf. \cite{T}).\\

\section{Proof of Theorem 1}
%First we need the following

%\begin{lem}
%For every $M\in \M$, and for every $a$ such that $e^a\in \sigma(S)$ we have\\
%1) $(Mf)_a\in L^2(\R),\:\forall f\in \Cc$.\\
%2) There exists $\mu_{(a)}\in L^\infty(\R)$ such that 
%$$\int_\R (Mf)(x)e^{ax}e^{-itx}dx=\mu_{(a)}(t)\int_\R f(x)e^{ax}e^{-itx}dx,\:a.e. $$
%i.e.\\
%$$\widehat{(Mf)_a}=\mu_{(a)} \widehat{(f)_a}.$$
%3) If $\overset{\circ}{I} \neq \emptyset $ then the function $\mu(z)=\mu_{(\Im z)}(\Re z)$ is holomorphic on $\overset{\circ}\Omega$.

For $\phi\in C_c(\R)$ denote by $M_\phi$ the operator of convolution by $\phi$ on $H$. We have
$$(M_\phi f)(x)=\int_\R f(x-y)\phi(y) dy,\:\forall f\in H.$$
It is clear that $M_\phi$ is a multiplier on $H$ for every $\phi\in \Cc$.\\
%Let $\A$ (resp. $\B$) be the closed algebra generated by operators $M_\phi$, for $\phi\in \Cc$ (resp. $S_x$, $x\in \R$) with respect to the topology of the operator norm. Denote by $\widehat{A}$ the set of characters of a commutatif algebra $A$.

%\end{lem}

In \cite{V4} we proved the following
\begin{thm}
%Let $E\subset L^1_{loc}(G)$ be a Banach space satisfying (H1), (H2) and (H3).
For every $M\in \M$, there exists a sequence $(\phi_n)_{n\in \N}\subset \Cc$ such that:\\
i) $ M=\lim_{n\to \infty} M_{\phi_n}$ with respect to the strong operator topology.\\
ii) We have $\|M_{\phi_n}\|\leq C\|M\|,$
%Every $M\in \M$ is the limit with respect to the strong operator topology of a net $(M_{\phi_\alpha})$ of operators of convolution with functions in $\CC$ such that 
%$$\|M_{\phi_\alpha}\|<C\|M \|,$$
where $C$ is a constant independent of $M$ and $n$. 
\end{thm}
The main difficulty to establish Theorem 1 is the proof of an estimate for $\widehat{\phi_n}(z)$ for $ z \in \Omega$ by the norm of $M_{\phi_n}.$

\begin{thm}
For every $\phi\in \Cc$ and every $\alpha \in \Omega$ we have
$$\Big|\int_\R \phi(x)e^{-i\alpha x} dx \Big|\leq \|M_\phi\|.$$
\end{thm}
Theorem 1 is deduced from Theorem 3 and Theorem 4 following exactly the same arguments as in Section 3 of \cite{V1} and Section 3 of \cite{V4}.
 The function $\mu_{(a)}$ introduced in Theorem 1 is obtained as the limit of $(\widehat{(\phi_n)_a})_{n\in \N}$
 with respect to the weak topology of $L^2(\R)$. The reader could consult \cite{V1} and \cite{V4} for more details.\\
 Here we give a proof of Theorem 4 by using the link between the spectrum of $S$ and the spectrum of the generator $A$ of the group $(S_t)_{t\in \R}$. \\

 {\bf Proof of Theorem 4.} 
 
 Let $\lambda\in \C$ be such that $e^{\lambda}\in \sigma(S)$. 
%First consider the case when $e^\lambda \in e^{\sigma(A)}.$ 
First we show that there exists a sequence $(n_k)_{k \in \N}$ of integers and a sequence $(f_{n_k})_{k\in \N}$ of functions of $H$ such that  
\begin{equation}\label{eq:1eq}
\|\Bigl(e^{t A} - e^{(\lambda + 2 \pi i n_k)t}\Bigr) f_{n_k}\| \longrightarrow 0,\: n_k \to \infty,\:\:\|f_{n_k}\|=1,\:\forall k\in \N.
\end{equation}
Let $A$ be the generator of the group $(S_t)_{t\in \R}$. 
We have to deal with two cases:\\
(i) $\lambda \in \sigma(A),$ \\
(ii) $\lambda \notin \sigma(A).$\\

In the case (i) we have $\lambda \in \sigma_p(A) \cup \sigma_c(A) \cup \sigma_r(A),$ where $\sigma_p(A)$ is the point spectrum, $\sigma_c(A)$ is the continuous spectrum and $\sigma_r(A)$ is the residual spectrum of $A$. If we have
$$\lambda \in \sigma_p(A) \cup \sigma_c(A),$$
it is easy to see that there exists a sequence $(f_m)_{m\in \N}\subset H$ such that 
$$\|(A-\lambda)f_m\|\underset{m \to +\infty}{\longrightarrow} 0,\:\|f_m\|=1,\:\forall m\in \N.$$
Then the equality 

$$(e^{At}-e^{\lambda t})f_m=\Big(\int_0^t e^{\lambda(t-s)} e^{As}ds\Big)(A-\lambda)f_m,$$
yields
$$\|(e^{At}-e^{\lambda t})f_m\|\underset{m \to +\infty}{\longrightarrow} 0, \:\forall t\in \R$$
and we obtain (\ref{eq:1eq}).
If $\lambda \notin \sigma_p(A)\cup \sigma_c(A)$, we have $\lambda \in \sigma_r(A)$ and 
$$\overline{Ran(A-\lambda I)}\neq H,$$
where $Ran(A-\lambda I)$ denotes the range of the operator $A-\lambda I$. 
 Therefore there exists $h\in D(A^*)$, $\|h\|=1,$ such that 
$$<f, (A^*-\overline{\lambda})h>=0,\: \forall f\in D(A).$$
This implies 
$(A^*-\overline{\lambda})h=0$ and we take $f=h$. 
Then
$$<(e^{At}-e^{\lambda t})f, f> =<f,(e^{A^*t}-e^{\overline{\lambda} t})f>$$
$$=\Big<f, \Big(\int_0^t e^{\overline{\lambda}( t-s)}e^{A^*s} ds\Big)(A^*-\overline{\lambda})f\Big>=0.$$
In this case we set $n_k=k$ and
$$f_{k}=f,\: \forall k\in \N$$
and we get again (\ref{eq:1eq}).\\
The case (ii) is more difficult since if $\lambda \notin \sigma(A)$, we have
$e^\lambda \in  \sigma(e^A)\setminus e^{\sigma(A)}.$
%Let $\w_0$ be the special weight equivalent to $\w_1$ (see the beginning of the section). Without lost of generality, we can consider $\{S_{t,0}\}$ and $A$ as operators on $L_{\w_0}^2(\R^2)$. 
%The spectrum of $S_{t,0}$ (resp. $A$) operating on $L_{\w_1}^2(\R^2)$ is the same as the spectrum of $S_{t,0}$ (resp. $A$) operating on $L_{\w_0}^2(\R^2)$. 
%Then, taking into account (\ref{eq:sg0}), we are in situation to apply in $L_{\w_0}^2(\R^2)$

Taking into account the results about the spectrum of a semi-group in Hilbert space \cite{H} satisfying the condition $(H4)$ (see also \cite{G} for the contraction semi-groups), we deduce that
there exists a sequence of integers
$n_k,$ such that $ |n_k|  \to \infty$ and
$$\|(A - (\lambda  + 2 \pi i n_k)I)^{-1}\| \geq k, \: \forall k \in \N.$$
Let $(g_{n_k})_{k\in \N}$ be a sequence such that 
$$\|g_{n_k}\| = 1,\: \Big\|\Bigl((A - (\lambda  + 2 \pi i n_k)I)^{-1}\Bigr)g_{n_k}\Big\| \geq k/2,\: \forall k \in \N.$$

We define 
$$f_{n_k} = \frac{  \Bigl((A- (\lambda  + 2 \pi i n_k)I)^{-1}\Bigr)g_{n_k}}{\|\Bigl((A - (\lambda  + 2 \pi i n_k)I)^{-1}\Bigr) g_{n_k}\|}.$$
Then we obtain
$$\Bigl(e^{t A}- e^{(\lambda + 2 \pi i n_k)t}\Bigr) f_{n_k} = \int_0^t e^{(\lambda + 2 \pi i n_k)(t-s)} e^{s A} ds \Bigl(A - (\lambda + 2 \pi i n_k)\Bigr) f_{n_k}$$
and for every $t$ we deduce
$$\lim_{k\to +\infty}\|\Bigl(e^{t A} - e^{(\lambda + 2 \pi i n_k)t}\Bigr) f_{n_k}\| =0.$$ 
Thus is established (\ref{eq:1eq}) for every $\lambda$ such that $e^{\lambda}\in \sigma(S)$.\\

Now consider
$$\hat{\phi}(-i\lambda) = \langle \int_\R \phi(t) \Bigl(e^{(\lambda + 2 \pi i n_k)t} - e^{t A}\Bigr) f_{n_k}, e^{ 2 \pi i n_k t} f_{n_k} \rangle dt
 + \langle \int_\R \phi(t) e^{t A} f_{n_k}, e^{ 2 \pi i n_k t} f_{n_k} \rangle dt $$
 $$=J_{n_k}+ \langle \int_\R \phi(t) e^{t A} f_{n_k}, e^{ 2 \pi i n_k t} f_{n_k} \rangle dt ,$$
where $J_{n_k} \to 0$ as $n_k \to \infty.$ 
On the other hand, we have
$$I_{n_k} = \langle \int_\R \phi(t) e^{t A} f_{n_k}, e^{ 2 \pi i n_k t} f_{n_k}\rangle dt = 
\langle \Bigl[ \int_\R \phi(t) e^{-2 \pi i n_k t} f_{n_k}(. - t)dt\Bigr], f_{n_k} \rangle$$
$$= \langle \int_\R \phi(. - y) e^{-2 \pi i n_k (.-y)} f_{n_k}(y) dy, f_{n_k}\rangle$$
$$ = \langle \Bigl(M_{\phi} ( f_{n_k} e^{ 2 \pi i n_k .})\Bigr), e^{2 \pi i n_k .} f_{n_k}\rangle$$
and $|I_{n_k}| \leq \|M_{\phi}\|.$ Consequently, we deduce that
$$|\hat{\phi}(-i\lambda)| \leq \|M_{\phi}\|.$$

Next a similar argument yields
\begin{equation}\label{eq:iv}
|\hat{\phi}(-i\lambda -a)| \leq \|M_{\phi}\|,\: \forall a \in \R.
\end{equation}

In fact, if for $t \in \R$ there exists a sequence $(h_n)_{n\in \N}\subset H$
 such that $(e^{tA} - e^{\lambda t})h_n \to 0$ as $n \to \infty$ with $\|h_n\| = 1$, we consider
$$< \int_\R(\phi(t) (e^{\lambda t} - e^{A t}))h_n, e^{-ia t} h_n> dt= \hat{\phi}(-i\lambda -a) - < \int_\R \phi(t) e^{i a t} e^{tA} h_n dt , h_n >.$$
The term on the left goes to 0 as $n \to \infty$, so it is sufficient to show that the second term on the right
is bounded by $\|M_{\phi}\|.$ We have
$$ \Big(\int_\R \phi(t) e^{i a t} e^{tA} h_n dt\Big)(x) = \int_\R \phi(t) e^{i a t} h_n(x - t) dt $$
$$ = \int_\R \phi(x- y) e^{i a (x-y)} h_n(y) dy = e^{ i a x} [M_{\phi}(e^{-ai .} h_n)](x),a.e.$$
and we obtain 
$$|\hat{\phi}(-i\lambda-a)|\leq \|M_{\phi}\|.$$\\
Next consider the second case when we have a sequence $(f_{n_k})_{k\in \N}$ with the properties above.
Multiplying by $e^{i (2 \pi n_k - a) t} f_{n_k}$, we obtain
$$\hat{\phi}(-i\lambda -a) = < \int_\R \phi(t) e^{t A} f_{n_k}, e^{ i (2 \pi  n_k - a) t} f_{n_k} > dt + I_{n_k},$$
where $I_{n_k} \to 0$ as $n_k \to \infty.$ 
To examine the integral on the right, we apply the same argument as above,  using the fact that  $(2 \pi n_k - a) \in \R.$ 
This completes the proof of (\ref{eq:iv}).
The property (\ref{eq:iv}) implies that if for some $\lambda_0\in\C$ we have 
$$|\hat{\phi}(\lambda_0)|\leq \|M_\phi\|,$$
then 
$$|\hat{\phi}(\lambda)|\leq \|M_\phi\|,\:\forall \lambda\in \C,\:s.t. \:\Im \lambda=\Im \lambda_0.$$

%This completes the proof of (1).

There exists $\alpha_0\in \sigma(S)$ such that $|\alpha_0|=\rho(S)$. Then we obtain that 
$$|\widehat{\phi}(z)|\leq \|M_\phi\|,$$
for every $z$ such that $\Im z =\ln \rho(S)$. 
In the same way there exists $\eta\in \sigma(S^{-1})$ such that $|\eta|=\rho(S^{-1})$ and 
$\alpha_1 = \frac{1}{\eta}\in \sigma(S)$. Then applying the above argument to $\alpha_1$, we get 
$$|\widehat{\phi}(z)|\leq \|M_\phi\|,$$
for every $z$ such that $\Im z =-\ln \rho(S^{-1})$. Since $\phi\in \Cc$ we have 
$$|\hat{\phi}(z)|\leq C\|\phi\|_{\infty}e^{k |\Im z|}\leq K \|\phi\|_{\infty},\:\:\forall z \in \Omega,$$
where $C > 0$, $k >0$ and $K >0$ are constants.
An application of the Phragmen-Lindel\"off theorem for the holomorphic function $\widehat{\phi}(z)$ yields
$$|\widehat{\phi}(\alpha)|\leq \|M_\phi\|$$
for all $\alpha\in \Omega$.
$\Box$\\
%\section{The spectrum of $S$}
Now we pass to the proof of Theorem 2. It is based on Theorem 1 combined with the arguments in \cite{V8} to cover our more general case. For the convenience of the reader we give the details.  \\

{\bf Proof of Theorem 2.}
Let $\alpha\in \C$ be such that $e^\alpha \notin \sigma(S)$. Then it is clear that $T=(S-e^{\alpha}I)^{-1}$ is a multiplier. Let $a\in \I$. Then there exists $\nu_{(a)}\in L^\infty(\R)$ such that
$$\widehat{(Tf)_a}=\nu_{(a)} \widehat{(f)_a},\: \forall f\in \Cc,\:a.e.$$
For $g\in \Cc$, the function $(S-e^\alpha I)g$ is also in $\Cc$. 
Replacing $f$ by $(S-e^\alpha I)g$, for $g\in \Cc$ we get
$$\widehat{(g)_a}(x)=\nu_{(a)}(x) \F\Big([(S-e^\alpha I)g]_a\Big)(x),\: \forall g\in \Cc,\:a.e.$$
and
$$\widehat{(g)_a}(x)=\nu_{(a)}(x)\widehat{g_a}(x)[e^{a-ix}-e^\alpha],\: \forall g\in \Cc,\:a.e.$$
Choosing a suitable $g\in \Cc$, we have
$$\nu_{(a)}(x)(e^{a-ix}-e^\alpha)=1, \:a.e.$$
On the other hand, $\nu_{(a)}\in L^\infty(\R)$. Thus we obtain that $\Re \alpha\neq a$ and we conclude that 

$$e^{a+ib}\in \sigma(S), \:\forall b\in \R.$$
Since $S$ is invertible, it is obvious that
$$\sigma(S)\subset \{z\in \C,\:\frac{1}{\rho(S^{-1})}\leq |z|\leq \rho(S)\},$$
Consequently, we obtain
$$\sigma(S)=\{z\in \C,\:\frac{1}{\rho(S^{-1})}\leq |z|\leq \rho(S)\}$$
and this completes the proof.
$\Box$

%On a vu que $\g$ est log-convexe et compact. Ces deux propriétés sont les seules contraintes géometriques sur $\g$.

%We conjecture that we have for every $M\in\M$, $\overline{\nu(\mathcal{U})}= \sigma(M).$

\end{document}